\documentclass[12pt,reqno]{amsart}
\newtheorem{theorem}{Theorem}

\theoremstyle{remark}

\def\ni{\noindent}

\def\be{\begin{equation}}
\def\ee{\end{equation}}

\def\bea{\begin{eqnarray}}
\def\eea{\end{eqnarray}}
\def\beas{\begin{eqnarray*}}
\def\eeas{\end{eqnarray*}}
\def\nn{\nonumber}

\def\oldfac#1{\,\hbox{\vrule width.5pt depth4pt
\hspace{-0.3mm}{\underline{\ $#1$\vphantom{gf}}}}}

\def\po#1#2{(#1)_#2}

\def\ep{\varepsilon}
\usepackage{amssymb}

\begin{document}

\newbox\Adr
\setbox\Adr\vbox{
\begin{center}
K. Srinivasa Rao$^{a,b}$, G. Vanden Berghe$^{a,c}$ and
C. Krattenthaler$^{d}$\\[.5cm]
$^a$Flemish Academic Center (VLAC), Royal Flemish Academy\\
of Belgium for Science and the Arts, Paleis der Academi\"{e}n,\\
Hertogsstraat 1, B-1000 Brussels, Belgium.\\[.2cm]
$^b$The Institute of Mathematical Sciences, C.I.T. Campus,\\
Taramani, Chennai - 600 113, India.\\[.2cm]
$^c$Universiteit Gent, Toegepaste Wiskunde en Informatica,\\
Krijgslaan 281-S9, B-9000 Gent, Belgium.\\[.2cm]
$^d$Institut Girard Desargues, Universit\'e Claude Berrnard Lyon-I,\\
21, avenue Claude Bernard,
F-69622 Villeurbanne, Cedex.\\[.5cm]
{\it E-mail}: {\tt rao@imsc.res.in, guido.vandenberghe@rug.ac.be,
kratt@euler.univ-lyon1.fr}
\end{center}
}

\newbox\Corr
\setbox\Corr\hbox{\footnotesize \noindent Corresponding author: G.
Vanden Berghe, Universiteit Gent, Toegepaste Wiskunde en
Informatica, Krijgslaan 281-S9, B-9000 Gent, Belgium}

\title[An entry of Ramanujan on Hypergeometric Series]
{\large An entry of Ramanujan on Hypergeometric Series in his
Notebooks\protect\footnote{\unhbox\Corr}}

\author[K. Srinivasa Rao, G. Vanden Berghe and C. Krattenthaler]
{\box\Adr}


\begin{abstract}
Example 7, after Entry 43, in Chapter XII of the first Notebook of
Srinivasa Ramanujan is proved and,
more generally, a summation theorem for
$_3F_2(a,a,x;1+a,1+a+N;1)$, where $N$ is a non-negative integer,
is derived.
\end{abstract}

\maketitle

\section{Introduction}

In the Notebooks of Ramanujan, identities for hypergeometric
series occupy a prominent part (see
\cite[Chapters X, XI]{Berndt}). In fact, Ramanujan discovered for
himself all the (now classical) summation theorems established by
Gau{\ss}, Chu-Vandermonde, Kummer, Pfaff--Saalsch\"utz, Dixon and
Dougall. For example, Dougall's summation theorem (of which all
the others mentioned are special cases), which was discovered by
Dougall \cite{Dougall} in 1907, was independently found by
Ramanujan during 1910--1912 (see Entry~1 in Chapter~XII of
Notebook~1 and in the corresponding Entry~1 in Chapter~X of
Notebook~2). We cannot assert an exact date, because there are no
dates anywhere in the Notebooks of Ramanujan. In any case,
Ramanujan rediscovered not only all that was known in Europe on
hypergeometric series at that time, but he also discovered several
new theorems, and, in particular, theorems on products of
hypergeometric series \cite{SlatAC}, as well as several types of
asymptotic
expansions.\\

On the other hand,
it is well known \cite{Berndt, Bradley} that Ramanujan did not
publish any of his results from his Chapters on hypergeometric
series that are to be found in his Notebooks. The Chapters on
hypergeometric series, in particular, Chapters~X and XI, in his
{\it second\/} Notebook, have been extensively studied by Hardy
\cite{Hardy} and Berndt \cite{Berndt}, since the Chapters of
Ramanujan's second Notebook are considered as {\it revised, enlarged\/}
versions of the Chapters in his {\it first\/} Notebook.
However, Example~7, after Entry~43, in Chapter XII of the first
Notebook did not find a place in the ``corresponding" examples
after Entry 10 in Chapter X of his second Notebook,
and is therefore not discussed in \cite{Berndt}.\\

The purpose of this article is to provide a proof for Example~7,
after Entry~43, in Chapter~XII, (XII, 43, Ex. 7), in the first
Notebook of Ramanujan, using well-known transformation and
summation theorems of hypergeometric series. In fact, we prove a
more general identity (which may have been the identity that
Ramanujan actually had, and from which he recorded the most
elegant specializations into his Notebook -- a conjecture
consistent with a wide spread belief that that was
Ramanujan's style).\\

Our paper is organized as follows: In Section 2, Entry (XII, 43,
Ex. 7) in the first Notebook of Ramanujan and some related entries
are presented in Ramanujan's notation, and subsequently translated
into current notation.  This leads us to the statement of a
summation theorem for a
${}_3F_2(x,\frac{1}{2},\frac{1}{2};\frac{3}{2},\frac{3}{2};1)$
hypergeometric series, for $\Re\, x<2$. In Section~3, the
statement is generalized to a summation theorem for the
${}_3F_2(a,a,x;1+a,1+a+N;1)$ series, where $N$ is a non-negative
integer, and it is proved. Finally, in Section~4, some remarks
regarding special cases of the theorem are made.

\section{Ramanujan's Entry}

In Ramanujan's notation, Example 7, after Entry 43, in
Chapter XII, of the first Notebook reads as follows:
\beas &&{\frac
\pi {\tan(\pi x)}}{\frac {\oldfac{2x}}
{\left(2^x\oldfac{x}\right)^2(1-2x)}} \left(\sum{\frac 1
{2x}}-{\frac 12}\sum {\frac1 x}+{\frac 1 {1-2x}}-{\frac\pi
2}\tan(\pi
x)\right)\\
&& \qquad = {\frac 1 {1^2}}+{\frac x {\oldfac{1}}}\cdot {\frac 1
{3^2}}+{\frac {x(x+1)} {\oldfac{2}}}\cdot {\frac 1 {5^2}}+\&c.
\qquad \qquad \qquad \qquad \qquad \kern20pt ({\rm XII}, 43, {\rm
Ex}.7)
 \eeas
Here, $\oldfac{x}$ is Ramanujan's notation for the {\it gamma
function} $\Gamma(x+1)$, which, for him, was a function over real
numbers $x$ (see \cite{Ramletter}).
We, of course, adopt the contemporary point of view and regard the
gamma function as a function over the complex numbers.\\

The factor on the left-hand side of (XII, 43, Ex.7): \be
{\frac\pi {\tan(\pi x)}}{\frac {\oldfac{2x}}
{\left(2^x\oldfac{x}\right)^2(1-2x)}} \label{2}\ee is the same
factor that appears on the left-hand side of (XII, 43, Ex.4)
which, in Ramanujan's first Notebook is
\[ \qquad \qquad {\frac \pi {\tan(\pi x)}}{\frac {\oldfac{2x}}
{\left(2^x\oldfac{x}\right)^2(1-2x)}} = 1+{\frac x
{\oldfac{1}}}\cdot {\frac 13}+{\frac {x(x+1)} { \oldfac{2}}}\cdot
{\frac 15}+ \&c. \qquad ({\rm XII}, 43, {\rm Ex}.4)\] This is
given in Ramanujan's second Notebook as:
\medskip

\noindent {\hfill $\displaystyle {\frac {\sqrt{\pi}\oldfac{n}} {
2\oldfac{n+\frac{1}{2}}}}=1-{\frac n {\oldfac{1}}}\cdot {\frac
13}+{\frac {n(n-1)} { \oldfac{2}}}\cdot {\frac 15}+ \&c. $ \hfill
({\rm X}, 10, {\rm Ex}.4)}

\medskip

 As pointed out by Berndt
\cite{Berndt}, the Entry (X, 10, Ex.4), or, (XII, 43, Ex.4), is
the special case of (X, 10), or, (XII, 43), where we do the
replacements $n\rightarrow \frac{1}{2}$ and
 $x\rightarrow n$. The factor (\ref{2}) which occurs on
the left-hand side of (XII, 43, Ex.4) can be shown to be equal to:
\be \frac {\sqrt{\pi}\,\Gamma(1-x)} {2x
\,\Gamma(\frac{3}{2}-x)}={\frac {\sqrt{\pi}\oldfac{-x}} {
2x\oldfac{-x+\frac{1}{2}}}} \;,\label{Faktor} \ee after using the
reflection formula: $\Gamma(z)\,\Gamma(1-z)=\pi / \sin(\pi z)$,
and the duplication formula
$\Gamma(2z)=2^{2z-1}\pi^{-1/2}\Gamma(z)\Gamma(z+\frac {1} {2})$
and some algebraic manipulations.\\

When written in the standard hypergeometric notation
\begin{equation*}  % \label{eq:hyper}
{}_r F_s\!\left[\begin{matrix} a_1,\dots,a_r\\
b_1,\dots,b_s\end{matrix}; z\right]=\sum _{k=0} ^{\infty}\frac
{\po{a_1}{k}\cdots\po{a_r}{k}} {k!\,\po{b_1}{k}\cdots\po{b_s}{k}}
z^k\ ,
\end{equation*}
where the Pochhammer symbol $(\alpha)_k$ is defined by $(\alpha)_k
=  \alpha(\alpha+1)(\alpha+2)\cdots (\alpha+k-1)$, $k>0$,
$(\alpha)_0  = 1$, the series on the right-hand side of (XII, 43,
Ex.4) is:
$$_2F_1\!\left[\begin{matrix} \frac{1}{2},x\\\frac{3}{2}\end{matrix};1
\right],$$ which by the Gau{\ss} summation theorem (see
\cite[(1.7.6); Appendix (III.3)]{SlatAC}): \bea\label{17} {} _{2}
F _{1} \!\left [ \begin{matrix} { a, b}\\ { c}\end{matrix} ;
{\displaystyle
   1}\right ]  = \frac { \Gamma (c)\,\Gamma( c-a-b)} {\Gamma( c-a)\,
\Gamma(
   c-b )} \quad \text{valid for }\Re(c-a-b)>0,
\eea becomes:
\be
{} _{2} F _{1} \!\left [ \begin{matrix} { \frac{1}{2}, x}\\
{ \frac{3}{2}}\end{matrix} ; {\displaystyle 1}\right]={\frac
{\sqrt{\pi}\,\Gamma(1-x)} {2\,\Gamma(3/2-x)}} ={\frac
{\sqrt{\pi}}2}{\frac {\oldfac{-x}}
  { \oldfac{-x+1/2}}}.
\label{corr}\ee
If we compare the right-hand side of this equation with (\ref{Faktor}),
which, as we outlined, is equal to the left-hand side (\ref{2}) of
(XII, 43, Ex.4),
it is clear that Ramanujan missed a
multiplicative factor $x$ on the left-hand side of (XII, 43,
Ex.4), while his Entry (X, 10, Ex.4) is correct.
As we shall see, the same applies to the Entry (XII, 43, Ex.7).
To be precise, the factor (\ref{2}) on the left-hand side of that entry
must also be replaced by the correct value (\ref{corr}).\\

The series on the right-hand side of (XII, 43, Ex.7)  is, in
hypergeometric notation:
$$ {} _{3} F _{2}
\!\left [ \begin{matrix} { x, \frac{1}{2}, \frac{1}{2}}\\
{ \frac{3}{2}, \frac{3}{2}}\end{matrix} ; {\displaystyle
1}\right].
$$
There are two factors on the left-hand side of (XII, 43, Ex.7).
Besides the factor (\ref{2}), the other factor in (XII, 43, Ex.7)
is: \be\sum{\frac 1 {2x}}-{\frac 12}\sum {\frac 1 x}+{\frac 1
{1-2x}}-{\frac \pi2}\tan(\pi x).\label{6}\ee Ramanujan used the
notation $\sum{\frac 1 x}$ to indicate the extension of the
function
$$1+\frac {1} {2}+\frac {1} {3}+\dots+\frac {1} {n},$$
representing the harmonic numbers, from positive integers $n$ to
real $x$. In other words, $\sum{\frac 1 x}$ is Ramanujan's
notation for the {\it digamma function}
$\psi(x):=\Gamma'(x)/\Gamma(x)$, the logarithmic derivative of the
gamma function, or, more precisely, for $\psi(x+1)+\gamma$, where
$\gamma$ is the Euler-Mascheroni constant, that is \be
\psi(x+1)=\sum \frac{1}{x}-\gamma.\label{9}\ee In fact, in
Ramanujan's very first research paper \cite{Ramanujan}, the
digamma function occurs as:
$$
\frac{d}{dn}\log \Gamma(n+1) = 1 +\frac{1}{2}+\frac{1}{3}+\cdots
+\frac{1}{n} - \gamma = \sum \frac{1}{n}-\gamma. %\label{8}
$$
Thus, the expression (\ref{6}) becomes: \be \psi(2x+1)-\frac {1}
{2}\psi(x+1)+ {\frac 1 {1-2x}}-{\frac \pi 2}\tan(\pi x)+\frac
\gamma 2\;.\label{7}\ee The digamma function satisfies the
recurrence relation \cite[1.7.1(8)]{Erdelyi}: \be
\psi(z+1)=\psi(z)+\frac{1}{z},\label{11}\ee the reflection formula
\cite[1.7.1(11)]{Erdelyi}: \be \psi(-z)=\psi(z+1) + \pi \cot(\pi z)
\label{10}\ee and the duplication formula \cite[1.7.1(12)]{Erdelyi}: \be
2\psi(2z)=\psi(z)+\psi\Big(z+\frac {1} {2}\Big)+2\log 2.
\label{dupli}\ee We remark that the duplication formula appears
implicitly in Ramanujan's first Notebook. Namely, a comparison of
Entry (XII, 43, Ex.6) with $n\rightarrow x-1$ in Entry (X, 10,
Ex.6) shows that Ramanujan obtained:
\begin{equation*}
\sum {\frac 1 {x+ {\frac 1 2}}}-\sum{\frac 1 x} =  2\sum {\frac 1
{2x}}-2\sum {\frac 1 x}+{\frac 1 x}-2\log 2,\end{equation*} which
by (\ref{9}) and (\ref{11}) is equivalent to (\ref{dupli}).
Furthermore, the digamma function has the special values \be
\psi\Big(\frac{1}{2}\Big)=-\gamma -2\log 2\qquad {\rm and} \qquad
\psi(1)=-\gamma.\label{12}\ee We now use the duplication formula
(\ref{dupli}) to convert (\ref{7}) into
$$\frac {1} {2}\psi\Big(x+\frac {1} {2}\Big)+\log 2+\frac {1} {1-2x}-
\frac {\pi} {2}\tan(\pi x)+\frac \gamma 2.$$
The recurrence relation (\ref{11}) implies that $\psi(x+\frac {1} {2})=
\psi(x-\frac {1} {2})+\frac {2} {2x-1}$, while we know from (\ref{12}) that
$\log 2=-\frac {1} {2}\psi(\frac {1} {2})-\frac {1} {2}\gamma$.
If this is substituted in the last equation, and if we then apply the
reflection formula (\ref{10}), we obtain
 \be
\frac {1} {2}\left(\psi\Big(\frac {3} {2}-x\Big)-
\psi\Big(\frac {1} {2}\Big)\right)\label{15}\ee
for the expression (\ref{6}). Therefore,
if we recall that the factor (\ref{2}) on the left-hand side of
Ramanujan's Entry (XII, 43, Ex.7) must be replaced by
(\ref{corr}), this entry can be rewritten in contemporary notation
as:
\begin{eqnarray}  {} _{3} F _{2} \!\left [
\begin{matrix} {x, \frac{1}{2}, \frac{1}{2}}\\ { \frac{3}{2},
\frac{3}{2}}\end{matrix} ; {\displaystyle 1}\right ]&=& {\frac
{\sqrt{\pi}} 4}{\frac {\Gamma(1-x)} {\Gamma({\frac 3 2}-x)}}
\left\{ \psi\Big({\frac 3 2}-x\Big)-\psi\Big({\frac 1 2}\Big)\right\}\nn\\
&=&{\frac 1 2}\left\{ \psi\Big({\frac 3 2}-x\Big)-\psi\Big({\frac 1
2}\Big)\right\}
\ {} _{2} F _{1} \!\left [ \begin{matrix} {x, \frac{1}{2}}\\
{ \frac{3}{2}}\end{matrix} ; {\displaystyle 1}\right ],\label{16}
\end{eqnarray} for $\Re \, x<2$.\\

In the next section we state and prove a theorem which is a
generalization of (\ref{16}), that is, of Ramanujan's Entry (XII,
43, Ex.7).

\section{The theorem}

\begin{theorem}
Let $N$ be a non-negative integer and $a$ be a complex number
which is not a negative integer. If $\Re \, x<N+2$, then
\begin{multline}
{} _{3} F _{2} \!\left [ \begin{matrix} { a, a , x}\\ { 1 + a , 1
+ a  + N}\end{matrix} ; {\displaystyle
1}\right ]\\
=\frac{a\,\Gamma(a + N+1)\,
       \Gamma(1 - x)\,
        }{N!\,\Gamma(a - x+1)}
\left( \psi(a - x+1)- \psi(a) -
         \psi(N+1) - \gamma
         \right)\\
- \frac{a\,\Gamma(a + N+1)\,
     \Gamma(1 - x)\,
     { }}{N!\,
     \Gamma(a - x+1)}
\sum _{k=1} ^ {N} \frac{({ \textstyle a}) _{k} \,({ \textstyle
-N}) _{k} }{k\cdot k! \,({ \textstyle
 a - x+1}) _{k}}.\label{18}
\end{multline}
\end{theorem}
\begin{proof}
To evaluate the ${}_3F_2(1)$ series on the left-hand
side of (\ref{18}), let us introduce a parameter $\ep$, and
consider the series
$${} _{3} F _{2} \!\left [ \begin{matrix} { a, a - \ep, x}\\ { 1 + a -
\ep, 1 + a - \ep + N}\end{matrix} ; {\displaystyle 1}\right ] .$$
First we apply the
(non-terminating) transformation formula (see \cite[Ex.~3.6, $q\to
1$, reversed]{GaRaAA}):
\begin{multline*}
{} _{3} F _{2} \!\left [ \begin{matrix} { a, b, c}\\ { d,
e}\end{matrix} ;
   {\displaystyle 1}\right ]
\\=
\frac {\Gamma( a - b)\,\Gamma( d)\,\Gamma( e)\,\Gamma(d + e -a - b
- c
 )} {\Gamma( a)\,\Gamma(  d-b)\,\Gamma(
     e -b )\,\Gamma(  d + e-a - c)}
\ {} _{3} F _{2} \!\left [ \begin{matrix} { b,  d-a,  e-a}\\ { 1 -
a + b,
     d + e-a - c}\end{matrix} ; {\displaystyle 1}\right ] \\+
 \frac { \Gamma( b-a)\, \Gamma( d)\, \Gamma( e)\,
 \Gamma( d + e-a - b - c  )}
{\Gamma( b)\, \Gamma(  d-a)\, \Gamma(
     e-a)\, \Gamma(  d + e-b - c)}
\ {} _{3} F _{2} \!\left [ \begin{matrix} { a,  d-b,  e-b}\\ {
 d + e-b - c, 1
     + a - b}\end{matrix} ; {\displaystyle 1}\right ].\label{18}
\end{multline*}
and obtain
\begin{multline*}
\frac{\Gamma({ \textstyle 1 + a - \ep}) \,
     \Gamma({ \textstyle \ep}) \,
     \Gamma({ \textstyle 1 + a - \ep + N}) \,
     \Gamma({ \textstyle 2 - \ep + N - x}) }{
\Gamma({ \textstyle a}) \,\Gamma({ \textstyle 1 + N}) \,
     \Gamma({ \textstyle 2 + a - 2 \ep + N - x}) }
\ {} _{2} F _{1} \!\left [ \begin{matrix} { a - \ep, 1 - \ep +
N}\\ { 2 + a - 2 \ep + N - x}\end{matrix} ;
{\displaystyle 1}\right ] \\
+
  \frac{\Gamma({ \textstyle 1 + a - \ep}) \,
     \Gamma({ \textstyle -\ep}) \,
     \Gamma({ \textstyle 1 + a - \ep + N}) \,
     \Gamma({ \textstyle 2 - \ep + N - x}) }{\Gamma({ \textstyle
1 - \ep}) \,\Gamma({ \textstyle a - \ep}) \,
     \Gamma({ \textstyle 1 - \ep + N}) \,
     \Gamma({ \textstyle 2 + a - \ep + N - x}) }\\
\times \ {} _{3} F _{2} \!\left [ \begin{matrix} { a, 1, 1 + N}\\
{ 1 + \ep, 2 + a - \ep + N - x}\end{matrix} ; {\displaystyle
1}\right ] .
\end{multline*}
%Christian: correction
Clearly, the convergence of the hypergeometric series on the
right-hand side will only be guaranteed if $\Re\,x<1$. Therefore. for
the moment we suppose $\Re\,x<1$.

To the $_3F_2$-series we apply the transformation formula (see
\cite[Ex.~7, p.~98]{BailAA}):
$$
{} _{3} F _{2} \!\left [ \begin{matrix} { a, b, c}\\ { d,
e}\end{matrix} ;
   {\displaystyle 1}\right ]  =
 \frac {\Gamma( e)\,\Gamma(d + e -a - b - c  )}
{\Gamma( e-a)\,\Gamma( d +
    e -b - c)}
 \ {} _{3} F _{2} \!\left [ \begin{matrix} { a,  d-b,  d-c}\\
 { d, d + e-b - c }\end{matrix} ; {\displaystyle 1}\right ] ,
$$
to get the expression
\begin{multline*}
  \frac{\Gamma({ \textstyle 1 + a - \ep}) \,
     \Gamma({ \textstyle \ep}) \,
     \Gamma({ \textstyle 1 + a - \ep + N}) \,
     \Gamma({ \textstyle 2 - \ep + N - x}) }{\Gamma({ \textstyle
1}) \,\Gamma({ \textstyle a}) \,\Gamma({ \textstyle 1 + N}) \,
     \Gamma({ \textstyle 2 + a - 2 \ep + N - x}) }
\ {} _{2} F _{1} \!\left [ \begin{matrix} { a - \ep, 1 - \ep +
N}\\ { 2 + a - 2 \ep + N - x}\end{matrix} ;
{\displaystyle 1}\right ]\\
 +
\frac{\Gamma({ \textstyle 1 + a - \ep}) \,
     \Gamma({ \textstyle -\ep}) \,
     \Gamma({ \textstyle 1 + a - \ep + N}) \,
     \Gamma({ \textstyle 1 - x}) }{\Gamma({ \textstyle 1 -
\ep}) \,\Gamma({ \textstyle a - \ep}) \,
     \Gamma({ \textstyle 1 - \ep + N}) \,
     \Gamma({ \textstyle 1 + a - x}) }
\ {} _{3} F _{2} \!\left [ \begin{matrix} { a, \ep, \ep - N}\\ { 1
+ \ep, 1 + a - x}\end{matrix} ; \ {\displaystyle 1}\right ].
\end{multline*}
The $_2F_1$-series is summed by means of the Gau{\ss} summation theorem
 (\ref{17}),
while the $_3F_2$-series is written as a sum over $k$, and
subsequently split into the ranges $k=0$, $k=1,2,\dots, N$, and
$k=N+1,N=2,\dots$. This yields the expression
\begin{multline*}
 \frac {1} {\ep}\bigg( \frac{\Gamma({ \textstyle 1 + a - \ep}) \,
     \Gamma({ \textstyle 1 + \ep}) \,
     \Gamma({ \textstyle 1 + a - \ep + N}) \,
     \Gamma({ \textstyle 1 - x}) }{
      \Gamma({ \textstyle a}) \,
     \Gamma({ \textstyle 1 + N}) \,
     \Gamma({ \textstyle 1 + a - \ep - x}) }  \\
- \frac{\Gamma({ \textstyle 1 + a - \ep}) \,
       \Gamma({ \textstyle 1 + a - \ep + N}) \,
       \Gamma({ \textstyle 1 - x}) }{
       \Gamma({ \textstyle a - \ep}) \,
       \Gamma({ \textstyle 1 - \ep + N}) \,
       \Gamma({ \textstyle 1 + a - x}) }\bigg)\\
-
  \frac{\Gamma({ \textstyle 1 + a - \ep}) \,
     \Gamma({ \textstyle 1 + a - \ep + N}) \,
     \Gamma({ \textstyle 1 - x}) \,
   }{
     \Gamma({ \textstyle a - \ep}) \,
     \Gamma({ \textstyle 1 - \ep + N}) \,
     \Gamma({ \textstyle 1 + a - x}) }
  \sum _{k=1} ^ {N} \frac{({ \textstyle a}) _{k} \,({ \textstyle \ep - N})
_{k} }{k! \,({ \textstyle k + \ep})
 \,({ \textstyle 1 + a - x}) _{k} }\\
 -
  \frac{\Gamma({ \textstyle 1 + a - \ep}) \,
     \Gamma({ \textstyle 1 + a - \ep + N}) \,
     \Gamma({ \textstyle 1 - x}) \,
  }{
     \Gamma({ \textstyle a - \ep}) \,
     \Gamma({ \textstyle 1 - \ep + N}) \,
     \Gamma({ \textstyle 1 + a - x}) }
   \sum _{k=N+1} ^ {\infty} \frac{({ \textstyle a}) _{k}
\,({ \textstyle \ep - N}) _{k} }{k!\,({ \textstyle k + \ep})
 \,({ \textstyle 1 + a - x}) _{k} }.
\end{multline*}
Now we perform the limit $\ep\to0$. Thus, our original
$_3F_2$-series becomes
$${} _{3} F _{2} \!\left [ \begin{matrix} { a, a , x}\\ { 1 + a ,
1 + a  + N}\end{matrix} ; {\displaystyle 1}\right ] .$$ On the
other hand, the four-line expression which we obtained for this
$_3F_2$-series simplifies significantly. The last term simply
vanishes because of the occurrence of the factor $({ \textstyle
\ep - N}) _{k}$, which is equal to
$(\ep-N)_N\,\ep\,(1+\ep)_{k-N-1}$ for $k\ge N+1$, which  is zero
for $N$ a non-negative integer. On the other hand, the limit of
$\ep\to0$ of the first two lines can be easily calculated by means
of l'H\^opital's rule to obtain the final result (\ref{18}).

As it stands, the assertion is only demonstrated for $\Re\,x<1$.
However, by analytic continuation, 
Equation~(\ref{18}) is true for any values
of $x$ for which the $_3F_2$-series on the left-hand side converges,
i.e., for $\Re\,x<N+1$.
\end{proof}

\eject

\newpage

\section{Some remarks}

\noindent The following observations can be made:
\begin{itemize}

\item[(i)] Clearly, the theorem (\ref{18}) reduces to Ramanujan's
entry (XII, 43, Ex.7), in our notation (\ref{16}), if $a=1/2$ and
$N=0$.\\

\item[(ii)] For $x=1$, in (\ref{16}), the $_3F_2(1)$ is a special
case of Dixon's theorem (see e.g.\ (III.8) of \cite{SlatAC}, for
$a=1,
b=c=1/2$), and it has the value ${\frac {\pi^2} 8}$.\\

\item[(iii)] For $x={\frac 3 2}$, the left-hand side of (\ref{16})
is: \beas _{3} F _{2} \!\left [ \begin{matrix} { \frac{3}{2},
\frac{1}{2},\frac{1}{2}}\\ { \frac{3}{2}, \frac{3}{2}}\end{matrix}
; {\displaystyle 1}\right ] = {}_2F_1 \left [ \begin{matrix} {
\frac{1}{2},\frac{1}{2}}\\ {\frac{3}{2}}\end{matrix} ;
{\displaystyle 1}\right ] = \frac{\pi}{2}, \eeas by the Gau{\ss}
summation theorem (\ref{17}), which is the result for the
right-hand side evaluated by l'H\^opital's rule.\\

\item[(iv)] For $x=-k$, a negative integer, in (\ref{16}), we get:
\be _{3} F _{2} \!\left [ \begin{matrix} { \frac{1}{2},
\frac{1}{2},-k}\\ { \frac{3}{2}, \frac{3}{2}}\end{matrix} ;
{\displaystyle 1}\right ]={\frac {\sqrt{\pi}} 2}{\frac
{\Gamma(k+1)} {\Gamma(3/2+k)}}\sum^{k+1}_{j=1}{\frac 1 {2j-1}}\;.
\label{20}\ee

\item[(v)] The $_3F_2(a,a,x;1+a,1+a+N;1)$ series
can be related to the 3-$j$ coefficient\break
$
\begin{pmatrix}
-\frac{x}{2} & -\frac{x}{2} & 0\\
a-\frac{x}{2} & -a+\frac{x}{2} & 0 \end{pmatrix} $ (cf.\
\cite{rao}) provided $x$ is a negative integer and $-x\leq a\leq
0.$ It also corresponds to the dual Hahn polynomial \cite{GARARR}
$S_n(0;a,b=1,c=1+N)$, for
$x=-n$.\\
\end{itemize}

Finally, as we already announced in Section~2,
it has to be noted that as in the case of Entry (XII, 43,
Ex.4), Ramanujan has missed a multiplicative factor $x$ on
the left-hand side of his Entry (XII, 43, Ex.7).\\

\ni {\bf Acknowledgments} : Two of the authors (KSR and GVB)
thank Prof.\ Dr.\ Niceas Schamp and the Koninklijke Vlaamse
Academie van Belgie voor Wetenschappen en Kunsten for excellent
hospitality and one of us (KSR) thanks Prof.\ Dr.\ Walter Van
Assche for an interesting discussion on a visit to the Katholieke
Universiteit
Leuven for a lecture on the Life and Work of Ramanujan.\\

\end{document}